\documentclass[11pt]{article}
\usepackage{amssymb, amsthm, amsmath, amscd}
\newtheorem{thm}{Theorem}[section]
\newtheorem{lem}[thm]{Lemma}
\newtheorem{prop}[thm]{Proposition}
\newtheorem{cor}[thm]{Corollary}
\newtheorem{NN}[thm]{}
\theoremstyle{definition}\newtheorem{df}[thm]{Definition}
\theoremstyle{definition}\newtheorem{rem}[thm]{Remark}
\theoremstyle{definition}\newtheorem{exm}[thm]{Example}

\newcommand{\N}{\mathbb{N}}
\newcommand{\Z}{\mathbb{Z}}
\newcommand{\Q}{\mathbb{Q}}
\newcommand{\R}{\mathbb{R}}

\newcommand{\T}{\mathbb{T}}

\newcommand{\hm}{homomorphism}
\newcommand{\dt}{\delta}
\newcommand{\ep}{\epsilon}
\newcommand{\andeqn}{\,\,\,{\rm and}\,\,\,}
\newcommand{\rforal}{\,\,\,{\rm for\,\,\,all}\,\,\,}
\newcommand{\CA}{$C^*$-algebra}

\newcommand{\beq}{\begin{eqnarray}}
\newcommand{\eneq}{\end{eqnarray}}
\newcommand{\tforal}{\,\,\,\text{for\,\,\,all}\,\,\,}
\newcommand{\tand}{\,\,\,\text{and}\,\,\,}

\title{The Range of Approximate Unitary Equivalence Classes of  Homomorphisms from AH-algebras}
\author{Huaxin Lin\\
 }
\date{}

\begin{document}

\maketitle

\begin{abstract}
Let $C$ be a unital AH-algebra and $A$ be a unital simple \CA\, with tracial rank zero.
It has been shown that two unital monomorphisms $\phi, \psi: C\to A$ are approximately unitarily equivalent
if and only if
$$
[\phi]=[\psi]\,\,\,{\rm in}\,\,\, KL(C,A)\andeqn \tau\circ \phi=\tau\circ \psi \tforal \tau\in T(A),
$$
where $T(A)$ is the tracial state space of $A.$ In this paper we
prove the following:  Given  $\kappa\in KL(C,A)$ with
$\kappa(K_0(C)_+\setminus \{0\})\subset K_0(A)_+\setminus \{0\}$
and with $\kappa([1_C])=[1_A]$ and a continuous affine map
$\lambda: T(A)\to T_{\mathtt{f}}(C)$ which is compatible with
$\kappa,$ where $T_{\mathtt{f}}(C)$ is the convex set of all
faithful tracial states,
 there exists a unital monomorphism
$\phi: C\to A$ such that
$$
[\phi]=\kappa\andeqn \tau\circ \phi(c)=\lambda(\tau)(c)
$$
for all $c\in C_{s.a.}$ and $\tau\in T(A).$ Denote by ${\rm Mon}_{au}^e(C,A)$ the set of approximate unitary equivalence classes of unital monomorphisms. We provide a bijective map
$$
\Lambda: {\rm Mon}_{au}^e (C,A)\to KLT(C,A)^{++},
$$
where $KLT(C,A)^{++}$ is the set of compatible pairs of elements
in $KL(C,A)^{++}$ and continuous affine maps from $T(A)$ to
$T_{\mathtt{f}}(C).$

Moreover, we realized  that there are compact metric spaces $X$,
unital simple AF-algebras $A$ and $\kappa\in KL(C(X), A)$ with
$\kappa(K_0(C(X))_+\setminus\{0\})\subset K_0(A)_+\setminus \{0\}$
for which there is no \hm\, $h: C(X)\to A$ so that $[h]=\kappa.$

\end{abstract}

\section {Introduction}

Recall that an  AH-algebra is a \CA\, which is an inductive limit of \CA s $C_n,$ where
$C_n=P_nM_{r(n)}(C(X_n))P_n$ for some finite CW complex $X_n$ and projections $P_n\in M_{r(n)}(C(X_n)).$
Note that every unital separable commutative \CA\, is an AH-algebra and every AF-algebra is an AH-algebra.
It was shown in \cite{Lncd} (see also Theorem 3.6 of \cite{Lnemb2}) that two unital monomorphisms
$\phi, \psi: C\to A,$ where $A$ is a unital simple \CA\, with tracial rank zero,
are approximately unitarily equivalent  if and only if
$$
[\phi]=[\psi]\andeqn \tau\circ \phi(c)=\tau\circ \psi(c)
$$
for all $c\in C_{s.a.}$ and $\tau\in T(A).$ This result plays a
role in the study of classification of amenable \CA s, or
otherwise known as the Elliott program. It also has applications
in the study of dynamical systems both classical and
non-commutative ones (see \cite{Lncd}). It is desirable to know
the range of the approximately unitary equivalence classes of
monomorphisms from a unital AH-algebra $C$ into a unital simple
\CA\, with tracial rank zero. For example, one may ask if given
any $\kappa\in KL(C,A)$ and any continuous affine map $\lambda:
T(A)\to T(C)$
 there exists a
monomorphism $\phi$ such that $[\phi]=\kappa$ and $\tau\circ
h(c)=\lambda(\tau)(c)$ for all $c\in C_{s.a.}$ and $\tau\in T(A).$

When $C$ is a finite CW complex, it was shown (see also a previous
result of L. Li \cite{Li}) in \cite{LnKT} that, for any $\kappa\in
KK(C,A)$ with $\kappa(K_0(C)_+\setminus\{0\})\subset
K_0(A)_+\setminus \{0\}$  and with $\kappa([1_C])=[1_A],$  there
exists a unital monomorphism $\phi: C\to A$ such that
$[\phi]=\kappa.$ It should be noted that both conditions that
$\kappa([1_C])=[1_A]$ and $\kappa(K_0(C)\setminus\{0\})\subset
K_0(A)\setminus \{0\}$ are necessary for the existence of such
$\phi.$ One of the earliest such results (concerning monomorphisms from $C(\T^2)$ into a unital simple AF-algebra) of this kind appeared in a paper of Elliott and Loring (\cite{EL} see also \cite{DL2}). It
was shown in \cite{LnKT} that the same result holds for the case
that $C$ is a unital simple AH-algebra which has real rank zero,
stable rank one and weakly unperforated $K_0(C).$ Therefore,
 it is natural to expect
that it holds for general unital AH-algebras.

 Let $C$ be the unitization of ${\cal K},$ the algebra of compact operators on $l^2.$
 Then it does not have a faithful tracial state. Consequently, it can not be
 embedded into
any unital  UHF-algebra, or any unital simple \CA\, which has  at
least one tracial state (It has been shown that a unital
AH-algebra $C$ can be embedded into a unital simple AF-algebra if
and only if $C$ admits a faithful tracial state --see
\cite{Lnemb2}) . This example at least suggests that for general
unital AH-algebras, the problem is slightly more complicated than
the first thought. Moreover, we note that to provide the range of
approximately unitary equivalence classes of unital monomorphisms
from $C,$ we also need to consider the map $\lambda: T(A)\to
T(C).$ Let $X$ be a compact metric space and let $C=C(X).$ Suppose
that $h: C\to A$ is a unital monomorphism and suppose that
$\tau\in T(A).$ Then $\tau\circ h$ induces a Borel probability
measure on $X.$ Suppose that $\kappa\in KL(C,A)$ is given. It is
clear that not every measure $\mu$ can be induced by those $h$ for
which $[h]=\kappa.$ Thus, we should consider a compatible pair
$(\kappa, \gamma)$ which gives a more complete information on
$K$-theory than either $\kappa$ or $\gamma$ alone.

The main result of this paper is to show that if $C$ is a unital
AH-algebra, $A$ is any unital simple \CA\, with tracial rank zero,
$\kappa\in KL(C,A)^{++}$ (see \ref{KL+} below) with
$\kappa([1_C])=[1_A]$ and $\lambda: T(A)\to T_{\mathtt{f}}(C), $
where $T_{\mathtt{f}}(C)$ is the convex  set of faithful tracial
states, which is a continuous affine  map and is compatible with
$\kappa,$ there is indeed  a unital monomorphism $\phi: C\to A$
such that
$$
[\phi]=\kappa\,\,\,{\rm in}\,\,\, KL(C,A)\andeqn  \phi_T=\lambda.
$$

We also show that the existence of $\lambda$ is essential to provide \hm s $\phi.$
In fact,
we find out  that there are compact metric spaces $X,$ unital simple AF-algebras $A$ and $\kappa\in KL(C(X),A)^{++}$
with $\kappa([1_C])=[1_A])$ for which there are  no $\lambda: T(A)\to T_{\mathtt{f}}(C(X))$ which is compatible with $\kappa.$   Moreover, we discovered  that there are no \hm\, $h: C\to A$ (not just monomorphisms) such that $[h]=\kappa.$
This further demonstrates that tracial information is an integral part of $K$-theoretical information.

\section{Notation}

\begin{NN}

{\rm Let $A$ be a unital \CA. Denote by $T(A)$ the tracial state space of $A.$ Denote by
${\rm Aff}(T(A))$ the space of all real affine continuous functions on $T(A).$ If $\tau\in T(A),$ we will also use $\tau$ for the tracial state $\tau\otimes {\rm Tr}$ on $M_k(A)$ for all integer $k\ge 1,$ where ${\rm Tr}$ is the standard trace on $M_k.$
If $a\in A_{s.a.},$ denote by $\check{a}$ a real affine function in ${\rm Aff}(T(A))$ defined by $\check{a}(\tau)=\tau(a)$
for all $\tau\in T(A).$

Let $C$ be another unital \CA. Suppose that $\gamma: {\rm Aff}(T(C))\to {\rm Aff}(T(A))$ is  a positive linear map.
We say it is unital if $\gamma(1_C)(\tau)=1.$ We say it is strictly positive, if $a\in {\rm Aff}(T(A))_+\setminus \{0\},$ then
$\gamma(a)(\tau)>0$ for all $\tau\in T(A).$

Suppose that $\phi: C\to A$ is a unital \hm. Denote by $h_T: T(A)\to T(C)$ the affine continuous map induced by $h,$ i.e.,
$$
h_T(\tau)(c)=\tau\circ h(c)\rforal c\in C.
$$
It also induces a positive linear map
$h_{\sharp}: {\rm Aff}(T(C))\to {\rm Aff}(T(A))$  defined by
$$
h_{\sharp}(\check{a})(\tau)=\tau\circ h(a)\tforal a\in C_{s.a}\andeqn \tau\in T(A),
$$
where $\check{a}(\tau)=\tau(a)$ for $\tau\in T(A).$

If $\lambda: T(A)\to T(C)$ is an affine continuous map, then it
gives a unital positive linear map $\lambda_{\sharp}: {\rm
Aff}(T(C))\to {\rm Aff}(T(A))$ by
$$
\lambda_{\sharp}(f)(\tau)=f(\lambda(\tau))\tforal f\in {\rm
Aff}(T(C))\andeqn \tforal \tau\in T(A).
$$
Conversely, a  unital positive linear map $\gamma: {\rm
Aff}(T(C))\to {\rm Aff}(T(A))$  gives an affine continuous map
$\gamma_T: T(A)\to T(C)$ by
$$
f(\gamma_T(\tau))=\gamma(f)(\tau)\tforal f\in {\rm
Aff}(T(C))\andeqn \tau\in T(C).
$$

Suppose that $A$ is a unital simple \CA. Then $\gamma$ is strictly
positive if and only if $\gamma_T$ maps $T(A)$ into
$T_{\mathtt{f}}(C).$

Denote by $\rho_A: K_0(A)\to {\rm Aff}(T(A))$ the positive \hm\, induced by
$\rho_A([p])(\tau)=\tau(p)$ for all projections
$p\in M_{\infty}(A)$ and $\tau\in T(A).$

Let $A$ and $C$ be two unital \CA s and let $\kappa_0: K_0(C)\to
K_0(A)$ be a unital positive \hm\, ($\kappa_0([1_C])=[1_A]$).
Suppose that $\lambda: T(A)\to T(C)$ is a continuous affine map.
We say that $\lambda$ is compatible with $\kappa_0,$ if
$\tau(\kappa([p]))=\lambda(\tau)(p)$ for all projections $p$ in
$M_{\infty}(A).$  Similarly, a unital positive linear map $
\gamma: {\rm Aff}(T(C))\to {\rm Aff}(T(A))$ is said to be
compatible with $\kappa_0,$ if
$\gamma(\check{p})(\tau)=\tau(\kappa([p])$ for all projections $p$
in $M_{\infty}(C).$ $\gamma$ is compatible with $\kappa_0$ if and
only if $\gamma_T$ is so.

Two projections in $A$ are equivalent if there exists a partial isometry $w\in A$ such that
$w^*w=p$ and $ww^*=q.$ }

\end{NN}

\begin{NN}\label{HomK}
{\rm Let $A$ be a unital \CA\, and let $C$ be a  separable \CA\,
which satisfies the universal coefficient theorem. By a result of
Dadarlat and Loring (\cite{DL1}),
\beq\label{N2-1}
KL(C,A)=Hom_{\Lambda}(\underline{K}(C), \underline{K}(A)),
\eneq
where, for any \CA\, $B,$
$$
\underline{K}(B)=\oplus_{i=0,1}K_i(B) \bigoplus_{n=2}^{\infty}\oplus_{i=0,1} K_i(B,\Z/n\Z).
$$
We will identify two objects in (\ref{N2-1}).
Denote by
$$
\underline{K}_{F,k}(C)=\bigoplus_{i=0,1}K_i(B) \bigoplus_{n|k}\oplus_{i=0,1} K_i(B,\Z/n\Z).
$$
If $K_i(C)$ is finitely generated ($i=0,1$), then there is $k_0\ge 1$ such that
$$
Hom_{\Lambda}(\underline{K}(C), \underline{K}(A))\cong 
Hom_{\Lambda}(F_{k_0}\underline{K}(C), F_{k_0}\underline{K}(A))
$$
(see \cite{DL1}). 
}

\end{NN}

\begin{df}\label{KL+}

{\rm
Denote by $KL(C,A)^{++}$ the set of those $\kappa\in Hom_{\Lambda}(\underline{K}(C), \underline{K}(A))$ such that
$$
\kappa(K_0(C)_+\setminus\{0\})\subset K_0(A)\setminus \{0\}.
$$
Denote by $KL_e(C,A)^{++}$ the set of those $\kappa\in
KL(C,A)^{++}$ such that $\kappa([1_C])=[1_A].$

}
\end{df}

\begin{df}\label{Dcomp}

{\rm Let $\kappa\in KL_e(C,A)^{++}$  and let  $\lambda: T(A)\to
T(C)$ be a continuous affine map. We say that $\lambda$ is
compatible with $\kappa$ if $\lambda$ is compatible with
$\kappa|_{K_0(C)}.$ Let $\gamma: {\rm Aff}(T(C))\to {\rm
Aff}(T(A))$ be a positive linear map. We say $\gamma$ is
compatible with $\kappa$ if  $\gamma$ is compatible with
$\kappa|_{K_0(C)},$ i.e., $\tau\circ
\kappa([p])=\gamma(\check{p})(\tau)$ for all projections $p\in
M_{\infty}(C).$ }
\end{df}

\begin{NN}\label{Meas}
{\rm Let $C=C(X)$ for some compact metric space $X.$ 
One has the following short exact sequence:
$$
0\to {\rm ker}\rho_C\to K_0(C)\to C(X, \Z)\to 0.
$$
It is then easy to see that, for every projection $p\in M_{\infty}(C),$ there is a projection $q\in C$ and an integer $n$ such that
$
\rho_A([p])=n\rho_A([q]).
$
It follows that if $C$ is a unital AH-algebra, then for every projection $p\in M_{\infty}(C),$ there is a projection
$q\in C$ and an integer $n\ge 1$ such that
$$
\rho_A([p])=n\rho_A([q]).
$$
Note also that in this case  ${\rm Aff}(T(C))=C_{s.a}.$
Therefore, in this note, instead of considering a unital positive linear maps $\gamma: {\rm Aff}(T(C))\to {\rm Aff}(T(A)),$ we may  consider unital positive linear maps $\gamma: C_{s.a}\to {\rm Aff}(T(A)).$ Moreover, $\gamma$ is compatible with some $\kappa\in KL(C,A)^{++},$ if
$\gamma(p)(\tau)=\tau(\kappa([p]))$ for all projections $p\in C$ and $\tau\in T(A).$

}
\end{NN}

\begin{NN}
{\rm Let $\phi, \psi: C\to A$ be two maps between \CA s. Let $\ep>0$ and ${\cal F}\subset C$ be a subset.
We write
$$
\phi\approx_{\ep} \psi\,\,\,{\rm on}\,\,\,{\cal F},
$$
if
$$
\|\phi(c)-\psi(c)\|<\ep\tforal c\in {\cal F}.
$$
}
\end{NN}

\begin{NN}
{\rm Let $L: C\to A$ be a linear map. Let $\dt>0$ and ${\cal G}\subset C$ be a (finite) subset. We say $L$ is $\dt$-${\cal G}$-multiplicative if
$$
\|L(ab)-L(a)L(b)\|<\dt\tforal a, b\in {\cal G}.
$$
}
\end{NN}

\begin{df}
{\rm  Let $A$ be a unital \CA. Denote by $U(A)$ the unitary group of $A.$ Let $B\subset A$ be another \CA\, and
$\phi: B\to A$ be a map. We write $\phi={\rm ad}\, u$ for some $u\in U(A)$ if
$\phi(b)=u^*bu$ for all $b\in B.$

Let $\phi, \psi: C\to A$ be two maps.  We say that $\phi$ and $\psi$ are approximately unitarily equivalent if there exists a sequence of unitaries $\{u_n\}\subset A$ such that
$$
\lim_{n\to\infty}{\rm ad}\, u_n\circ \phi(c)= \psi(c)\tforal c\in C.
$$

}
\end{df}

\section{Approximate unitary equivalence}

We begin with the following theorem

\begin{thm}{\rm (Theorem 3.6 of \cite{Lnemb2} and see also Theorem 3.4 of \cite{Lncd}) } \label{Om}
Let $C$ be a unital AH-algebra and let $A$ be a unital simple \CA\, with tracial rank zero.
Suppose that $\phi, \psi: C\to A$ are two unital monomorphisms.
Then there exists a sequence of unitaries $\{u_n\}\subset A$ such that
$$
\lim_{n\to\infty}{\rm ad}\, u_n\circ \psi(c)=\phi(c)\tforal c\in C,
$$
 if and only if
$$
[\phi]=[\psi]\,\,\,{\rm in}\,\,\,KL(C,A)\tand \tau\circ
\phi=\tau\circ \psi \tforal \tau\in T(A).
$$
\end{thm}

We need the following variation of results in \cite{Lncd}.

\begin{thm}\label{1Tcd}
Let $C$ be a unital $AH$-algebra, let $A$ be a unital simple \CA\, with $TR(A)=0$ and let  $\gamma: C_{s.a.}\to
{\rm Aff}(T(A))$ be a unital strictly positive linear map.

For any $\ep>0$ and any finite subset ${\cal F}\subset C,$ there exist $\eta>0,$  $\dt>0,$ a finite subset
${\cal G}\subset C,$ a finite subset ${\cal H}\subset C_{s.a.}$ and a finite subset ${\cal P}\subset \underline{K}(C)$ satisfying the following:

Suppose that  $L_1, L_2: C\to A$  are two  unital completely positive linear maps which are $\dt$-${\cal G}$-multiplicative  such that
 \beq\label{1Tcd-1}
[L_1]|_{\cal P}&=&[L_2]|_{\cal P},\\
|\tau\circ L_i(g)-\gamma(g)(\tau)|&<&\eta\tforal g\in {\cal H},\,\,\,i=1,2.
\eneq
Then there is a unitary $u\in A$ such that
\beq\label{1Tcd-2}
{\rm ad}\, u\circ L_2\approx_{\ep} L_2\,\,\,{\rm on}\,\,\,{\cal F}.
\eneq
\end{thm}

\begin{proof}
Write $C={\overline{\cup_{n=1}^{\infty}C_n}},$ where $C_n=P_nM_{r(n)}(C(X_n))P_n,$ where $X$ is a compact 
subset
of a finite CW complex and where $P_n\in M_{r(n)}(C(X_n))$ is a projection.
Let $\ep>0$ and a finite subset ${\cal F}\subset C$ be fixed.  Without loss of generality, we may assume
that ${\cal F}\subset C_1.$ Let $\eta_0>0$ such that
$$
|f(x)-f(x')|<\ep/8\tforal f\in {\cal F},
$$
if ${\rm dist}(x,x')<\eta_0.$
Let $\{x_1,x_2,...,x_m\}\subset X$ be $\eta_0/2$-dense in $X.$ Suppose that
$O_i\cap O_j=\emptyset$ if $i\not=j,$ where
$$
O_j=\{x\in X: {\rm dist}(x,x_j)<\eta_0/2s\},\,\,\,j=1,2,...,m
$$
for some integer $s\ge 1.$

Choose non-zero element $g_j\in (C_1)_{s.a}$ such that $0\le g_j\le 1$ whose support lies in $O_j,$
$j=1,2,...,m.$ Note such $g_j$ exists (by taking those in the center for example).
Choose
$$
\sigma_0=\min\{\inf\{\gamma(g_j)(\tau): \tau\in T(A)\}: 1\le j\le m\}.
$$
Since $\gamma$ is strictly positive, $\sigma_0>0.$  Set $\sigma=\min\{\sigma_0/2, 1/2s\}.$
Then, by Corollary 4.8 of \cite{Lncd}, such $\dt>0,$ $\eta>0,$ ${\cal G},$ ${\cal H}$ and ${\cal P}$ exists.

\end{proof}

\begin{lem}\label{Lcd}
Let $X$ be a compact metric space, let $A$ be a unital simple
\CA\, with $TR(A)=0$  and let $\gamma:C(X)_{s.a.}\to {\rm Aff}(T(A))$ be
a unital strictly positive linear map.

Then, for any $\ep>0$ and any ${\cal F}\subset C(X),$ there exists
$\dt>0,$ a finite subset ${\cal G}\subset C(X)_{s.a.},$ a set
 $S_1, S_2,...,S_n$ of mutually disjoint clopen subsets
 with $\cup_{i=1}^n S_i=X,$ satisfying the following:

 For any two unital \hm s $\phi_1, \phi_2: C(X)\to pAp$ with
 finite dimensional range for some projection $p\in A$  with $\tau(1-p)<\dt$
 such that
\beq\label{lcd-1}
[\phi_1(\chi_{S_i})]&=&[\phi_2(\chi_{S_i})]\,\,\,{\rm in}\,\,\,
K_0(A),\,\,\,i=1,2,...,n,\\
|\tau\circ \phi_1(g)-\gamma(g)(\tau)|&<&\dt\tand\\
|\tau\circ \phi_2(g)-\gamma(g)(\tau)|&<&\dt
\eneq
for all $g\in {\cal G}$  and for all $\tau\in T(A),$ there exist a
unitary $u\in U(pAp)$ such that
\beq\label{lcd-2}
{\rm ad}\, u\circ \phi_1\approx_{\ep} \phi_2\,\,\,{\rm on}\,\,\,
{\cal F}.
\eneq

\end{lem}

\begin{proof}
This follows from \ref{1Tcd} immediately. There is a sequence of
finite CW complex $X_n$ such that
$C(X)=\lim_{n\to\infty}(C(X_n),h_n),$ where each $h_n$ is a unital
\hm. Fix $\ep>0$ and a finite subset ${\cal F}\subset C(X).$
Without loss of generality, we may assume that ${\cal F}\subset
h_K({\cal F}_K)$ for some integer $K\ge 1$ and a finite subset
${\cal F}_K.$

Given any finite subset ${\cal P}\subset \underline{K}(C(X)),$ one
obtains a finite subset ${\cal Q}_k\subset \underline{K}(C(X_k))$
such that $[h_k]({\cal Q}_k)={\cal P}$ for some $k\ge 1.$ Let
$p_1,p_2,...,p_n$ be mutually orthogonal projections corresponding
to the connected components of $X_k.$    To simplify notation,
without loss of generality, we may  assume that $k=K.$

There are mutually disjoint clopen sets $S_1, S_2,...,S_n$ of $X$
with $\cup_{i=1}^n S_i=X$ such that $h_k(p_i)=\chi_{S_i},$
$i=1,2,...,n.$ Since $\phi$ and $\psi$ are \hm s with finite
dimensional range, if
$$
[\phi(\chi_{S_i})]=[\psi(\chi_{S_i})]\,\,\,{\rm in}\,\,\, K_0(A),
$$
then
$$
[\phi\circ h_k]=[\psi\circ h_k]\,\,\,{\rm in}\,\,\, KL(C(X_k),A).
$$
This, in particular, implies that
$$
[\phi]|_{\cal P}=[\psi]|_{\cal P}.
$$

This above argument shows that the lemma follows from \ref{1Tcd}.

\end{proof}

\begin{df}\label{D1}
Let $X$ be a compact metric space which is a compact subset of
some finite CW complex $Y.$ Then there exists a decreasing
sequence of finite CW complexes $X_n\subset Y$ such that
$$
X\subset X_n\andeqn \lim_{n\to\infty}\text{dist}(X_n, X)=0.
$$
Denote by $s_{m,n}: C(X_m)\to C(X_n)$ (for $n>m$) and $s_n:
C(X_n)\to C(X)$ be the surjective \hm s induced by the inclusion
$X_{n+1}\subset X_n$ and $X\subset X_n,$ respectively.
\end{df}

\begin{lem}\label{KLL1}
Let $Y$ be a finite CW complex and $X\subset Y$ be a compact
subset. For any $\ep>0,$ any finite subset ${\cal F}\subset C(X),$
there exists a finite subset ${\cal P}\subset
\underline{K}(C(X)),$  an integer $k\ge 1$ and an integer $N\ge 1$
satisfying the following:

For any unital homomorphisms $\phi, \psi: C(X_m)\to A$ ($m\ge k$)
for any unital simple \CA\, with $TR(A)=0$  for which
$$
[\phi]|_{\cal Q}=[\psi]|_{\cal Q}\,\,\,{\rm in}\,\,\, KL(C(X_m),
A),
$$
where $Q\subset \underline{K}(C(X_m))$ is a finite subset such
that $[s_m]({\cal Q})={\cal P},$ then there exists a unitary $U\in
M_{N+1}(A)$ such that
$$
\text{ad}\, U\circ (\phi\oplus \Phi\circ s_m)\approx_{\ep}
(\psi\oplus \Phi\circ s_m)\,\,\,\text{on}\,\,\, s_m^{-1}({\cal F}),
$$
 where $\Phi: C(X)\to M_N(A)$ is
defined by
\beq\label{KLL1-0}
\Phi(f)={\rm diag}(f(x_1), f(x_2),...,f(x_N))\tforal f\in
C(X_1),
\eneq
where $\{x_1, x_2,...,x_N\}$ is a finite subset of $X.$
\end{lem}

\begin{proof}
Assume that the lemma were false. Then  there would be  a positive
number $\ep_0>0,$ a finite subset ${\cal F}_0\subset C(X),$
an increasing sequence of finite subsets $\{{\cal P}_n\}\subset
\underline{K}(C(X))$ with $\cup_n{\cal P}_n=\underline{K}(C(X)),$
a sequence of unital \CA s, two subsequences $\{R(n)\},\,\{k(n)\}$
of $\N$ and two sequences monomorphisms $\phi_n, \psi_n:
C(X_{k(n)})\to A_n$ such that
\beq\label{KLL1-1}
&&[\phi_n]|_{{\cal Q}_n}=[\psi_n]|_{{\cal Q}_n}\,\,\,
\text{in}\,\,\, KK(C(X_{k(n)}), A_n)\andeqn\\\label{N-0}
&&\hspace{-0.8in}\limsup_n\{\inf\{\max\{\|u_n^*( \phi_n\oplus
\Phi_n\circ s_n)(f))u_n-(\phi\oplus \Phi_n\circ s_n)(f)\|: f\in s_m^{-1}({\cal F})\}\}\}\ge
\ep_0,
\eneq
where infimum is taken among all possible $\Phi_n:C(X)\to
M_{R(n)}(A_n)$ with the form described above and among all
possible unitaries $\{u_n\}\subset U(M_{R(n)+1}(A)),$ and where 
${\cal Q}_n\subset \underline{K}(C(X_{k(n)}))$ is a finite subset such that
$[s_{k(n)}]({\cal Q}_n)={\cal P}_n.$
Since $K_i(C(X_n)$ is finitely generated, by passing to a subsequence, if necessary, 
without loss of generality, we may assume (see also the end of \ref{HomK}) that 
\beq\label{N-00}
[\phi_{n+1}\circ s_{k(n), k(n+1)}]=[\psi_{n+1}\circ s_{k(n), k(n+1)}]\,\,\,{\rm in}\,\,\, KL(C(X_{k(n)}), A),\,\,\,n=1,2,....
\eneq
Let $\phi_n^{(m)}=\phi_m,$ if $n\le m,$ $\phi_n^{(m)}=\phi_n\circ
s_{m,n},$ $\psi_n^{(m)}=\psi_m,$ if $n\le m$ and
$\psi_n^{(m)}=\psi_n\circ s_{m,n},$ $n=1,2,....$ Denote by
$H_1^{(m)}, H_2^{(m)} : C(X_{k(m)})\to \prod_{n}A_n$ by
$H_1^{(m)}(f)=\{\phi_n^{(m)}\}$ and
$H_2^{(m)}(f)=\{\psi_n^{(m)}\}.$ Let $\pi: \prod_nA_n\to
\prod_nA_n/\bigoplus_nA_n$ be the quotient map. Then
 $\pi\circ H_1^{(m)}$ and $ \pi\circ H_2^{(m)}$ both have
spectrum $X.$  Moreover, for each $i,$ all $\pi\circ H_i^{(m)}$ gives the same 
\hm\, $F_i: C(X)\to 
\prod_n A_n/\bigoplus_n A_n,$ $i=1,2.$

Since $TR(A_n)=0,$ $A_n$ has real rank zero, stable rank one, weakly unperforated $K_0(A_n),$
by Corollary 2.1 of \cite{GL2} and (\ref{N-00})
$$
[H_1^{(m+1)}\circ s_{k(m), k(m+1)}]=[H_2^{(m+1)}\circ s_{k(m), k(m+1)}]\,\,\,{\rm in}\,\,\, KL(C(X_{k(m)}), \prod_nA_n)
$$

It follows from Corollary 2.1 of \cite{GL2} again  that
$$
[F_1]=[F_2]\,\,\,{\rm in}\,\,\, KL(C,\prod_nA_n/\bigoplus_nA_n).
$$
It then follows from  Theorem 1.1 and the Remark 1.1 of \cite{GL2} that there is an integer $N\ge 1$ and a
unitary $W\in U(M_{N+1}(\prod_n A_n/\bigoplus_n A_n))$ such that
\beq\label{N-2}
{\rm ad}\,W\circ (F_2\oplus H_0)\approx_{\ep_0/2} (F_1\oplus
H_0)\,\,\,{\rm on}\,\,\, {\cal F}_0,
\eneq
where $H_0: C(X)\to M_N(\prod_n A_n/\bigoplus_n A_n)$ is defined by
$H_0(f)=\sum_{i=1}^N f(x_i)E_i$ for all $f\in C(X),$ $x_i\in X$
and $E_i={\rm diag}(\overbrace{0,...,0}^{i-1}, 1, 0,...,0),$
$i=1,2,...,N.$

There is a unitary $\{W_n\}\in U(\prod_n A_n)$ such that
$\pi(\{W_n\})=W.$ Then, for some sufficiently large $n,$
\beq\label{N-3}
W_n^*{\rm diag}(\phi_n(f), f(x_1), f(x_2),...,f(x_N))W_n
\approx_{\ep_0}(\psi_n(f), f(x_1), f(x_2),...,f(x_N))
\eneq
on ${\cal F}_0.$ This contradicts (\ref{N-0}).

\end{proof}

\begin{rem}\label{RM1}
{\rm
There exists a positive number $\eta>0$ and integer $N_1>0$ which depend only on $\ep$ and ${\cal F}$ such that
$\{x_1,x_2,...,x_N\}$ and an integer $N$ can be replaced by any $\eta$-dense finite subset $\{\xi_1,\xi_2,...,\xi_{N_1}\}$
and integer $N_1.$

From the proof, we also know that the assumption that $A$ has
tracial rank zero can be replaced by much weaker conditions (see
Corollary 2.1 of \cite{GL2}). The main difference of \ref{KLL1}
and results in \cite{GL2} is that \hm s $\phi$ and $\psi$ are not
assumed to be from $C(X)$ to $A.$

}

\end{rem}

\section{Monomorphisms from $C(X)$}

\begin{lem}\label{KLL2}
Let $X$ be a finite CW complex and let $A$ be a unital simple
\CA\, with real rank zero, stable rank one and weakly unperforated
$K_0(A).$  Let $e_1,e_2,...,e_m\in C(X)$ be mutually orthogonal
projections corresponding to connected components of $X.$

Suppose that $\kappa\in KK(C(X), A)^{++}$ with
$\kappa([1_{C(X)}])=[1_A].$ Then, for any projection $p\in A$ and
any unital \hm\, $\phi_0: C(X)\to (1-p)A(1-p)$ with finite
dimensional range such that $\phi_0([e_i])< \kappa([e_i]),$
$i=1,2,...,m.$  Then there exists a unital monomorphism $\phi_1:
C(X)\to pAp$ such that
\beq\label{KLL2-1}
[\phi_1+\phi_0]=\kappa\,\,\,\text{in}\,\,\, KK(C(X), A).
\eneq
\end{lem}

\begin{proof}
Since $\sum_{i=1}^m\kappa([e_i])=[1_A]$ and $A$ has stable rank
one, there are mutually orthogonal projections $p_1,p_2,...,p_m\in
A$ such that
\beq\label{L2-2}
\sum_{i=1}^m p_i=1_A\andeqn [p_i]=\kappa([e_i]),\,\,\,i=1,2,...,m
\eneq
From this it is clear that we may reduce the general case to the
case that $X$ is connected.

So now we assume that $X$ is connected. Then it is easy to see
that
$$
\kappa-[\phi_0]\in KK(C(X), A)^{++}
$$
and $(\kappa-[\phi_0])([1_{C(X)}])=p.$ It follows from Theorem 4.7
of \cite{LnKT} that there is a monomorphism $\phi_1: C(X)\to pAp$
such that
$$
[\phi]=\kappa-[\phi_0].
$$

\end{proof}

\begin{lem}\label{ML0}
Let $X$ a compact metric space and let $A$ be a unital simple
\CA\, with tracial rank zero. Suppose that $\gamma: C(X)_{s,a} \to
{\rm Aff}(T(A))$ is a unital strictly positive linear map. Let
$S_1,S_2,...,S_n$ be a set of mutually disjoint clopen
subsets of $X$ with $\cup_{i=1}^n S_i=X.$  Then for any $\dt>0$
and any finite subset ${\cal G}\subset C(X)_{s.a},$ there exists a
projection $p\in A$ with $p\not=1_A$ and a unital \hm\, $h:
C(X)\to pAp$ with finite dimensional range such that
\beq\label{M0-1}
|\tau\circ h(g)-\gamma(g)(\tau)|&<&\dt\tforal g\in {\cal G}\tand
\tau\in T(A), \tand\\\label{M0-1+} \tau\circ
h(\chi_{S_i})&<&\gamma(\xi_{S_i})(\tau)\tforal \tau\in T(A),
\eneq
$i=1,2,...,n.$
\end{lem}

\begin{proof}
Put
$$
d=\min\{\dt, \min\{\inf\{\gamma(\chi_{S_i})(\tau):\tau\in
T(A)\}:1\le i\le n\}\}.
$$
Since $\gamma$ is strictly positive, $d>0.$

 Let ${\cal G}_0={\cal G}\cup \{\chi_{S_1}, \chi_{S_2},...,\chi_{S_n}\}.$
  It follows from 4.3  of \cite{Lnhomp}
that there is a unital \hm\, $h_0: C(X)\to A$ with finite dimensional range such that
\beq\label{M0-2}
|\tau\circ h(g)-\gamma(g)(\tau)|<d/8n \tforal g\in {\cal G}_0
\eneq
and for all $\tau\in T(A).$ In particular,
\beq\label{M0-3}
|\tau\circ h(\chi_{S_i})-\gamma(\chi_{S_i})(\tau)|<d/8n\tforal
\tau\in T(A)
\eneq
$i=1,2,...,n.$


Since $\rho_A(K_0(A))$ is dense in ${\rm Aff}(T(A)),$ there
exists a projection $p_0\in A$ such that
\beq\label{M0-4}
d/2n<\tau(p_0)<d/n \tforal \tau\in T(A).
\eneq

Note that $\tau(p_0)<\gamma(\chi_{S_i})(\tau)$ for all $\tau\in
T(A),$ $i=1,2,...,n.$ Moreover, by (\ref{M0-3}), 
\beq\label{M0-4+}
\tau\circ h(\chi_{S_i})> \gamma(\xi_{S_i})(\tau)-d/8n\ge
d-d/8n>\tau(p_0).
\eneq
for all $\tau\in T(A).$

Write $h_0(f)=\sum_{k=1}^m f(x_k)e_k$ for all $f\in C(X),$ where
$x_k\in X$ and $\{e_1,e_2,...,e_k\}$ is a set of mutually
orthogonal projections with $\sum_{k=1}^m e_k=1_A.$

Note that
$$
h_0(\xi_{S_j})=\sum_{x_k\in S_j}e_k.
$$
Therefore (by (\ref{M0-4+}))
\beq\label{M0-5}
[p_0]\le [\sum_{x_k\in S_j}e_k].
\eneq
By Zhang's Riesz interpolation property (see \cite{Z1} ), there are
projections $e_k'\le e_k$ such that
$$
[p_0]=[\sum_{k\in S_j}e_k'].
$$
By Zhang's half projection theorem (see Theorem 1.1 of \cite{Z2}),
for each $k,$ there is a projection $e_k''\le e_k'$ such that
\beq\label{M0-6}
[e_k'']+[e_k'']\ge [e_k'].
\eneq
Thus
\beq\label{M0-7}
2[\sum_{\chi_k\in S_i}e_k'']\ge [p_0],\,\,\,i=1,2,...,n.
\eneq
Therefore (by (\ref{M0-4}) and (\ref{M0-3}))
\beq\label{M0-8}
\tau(\sum_{x_k\in S_i}(e_k-e_k''))&<&\tau\circ
h_0(\chi_{S_i})-(1/2)\tau(p_0)\\
&<& \tau\circ h(\chi_{S_i})-d/4n\\\label{M0-8+}
&<& \gamma(\chi_{S_i})(\tau)-d/8n\tforal \tau\in T(A).
\eneq

Let $p=\sum_{k=1}^m (e_k-e_k'').$ Then clearly that $p\not=1.$
Moreover,
$$
\tau(1-p)<d/4\tforal \tau\in T(A).
$$

Define $h(f)=\sum_{k=1}^m f(x_k) (e_k-e_k'')$ for all $f\in C(X).$
Then
\beq\label{M0-9}
|\tau\circ h(f)-\tau\circ h_0(f)|<\tau(\sum_{k=1}^m
e_k'')=\tau(1-p)<d/4<\dt
\eneq
for all $\tau\in T(A).$

Then, by (\ref{M0-8+}),
\beq\label{M0-10}
\tau\circ h(\chi_{S_i})<\gamma(\chi_{S_i})(\tau)\tforal \tau\in
T(A).
\eneq

\end{proof}

\begin{lem}\label{ML1}
Let $X$ a compact metric space  and let $A$ be a unital simple
\CA\, with tracial rank zero. Suppose that $\gamma: C(X)_{s,a} \to
{\rm Aff}(T(A))$ is a unital strictly positive linear map. Let
$S_1,S_2,...,S_n$ be a set of mutually disjoint clopen
subsets of $X$ with $\cup_{i=1}^n S_i=X.$  Then for any $\dt>0,$
$\eta>0,$  for any integer $N$ and any $\eta$-dense subset
$\{x_1,x_2,...,x_N\}$ of $X$ and any finite subset ${\cal
G}\subset C(X)_{s.a},$ there exists a projection $p\in A$ with
$p\not=1_A$ and a unital \hm\, $h: C(X)\to pAp$ with finite
dimensional range such that
\beq\label{M1-1}
|\tau\circ h(g)-\gamma(g)(\tau)|&<&\dt\tforal g\in {\cal G}\tand
\tau\in T(A), \tand\\\label{M1-1+} \tau\circ
h(\chi_{S_i})&<&\gamma(\chi_{S_i})(\tau)\tforal \tau\in T(A),
\eneq
$i=1,2,...,n,$
\beq\label{M1-2}
h(f)=\sum_{i=1}^N f(x_i)e_i\oplus h_1(f)\tforal f\in C(X),
\eneq
where $h_1: C(X)\to (1-\sum_{i=1}^Ne_i)A(1-\sum_{i=1}^N e_i)$ is a
unital \hm\, with finite dimensional range and
$\{e_1,e_2,...,e_N\}$ is a set of mutually orthogonal projections
such that $[e_i]=[e_1]\ge [1-p],$ $i=1,2,...,N.$

\end{lem}

\begin{proof}
  Let $N\ge 1$ and let
$\eta$-dense subset $\{x_1, x_2,...,x_N\}$ of $X$ be  given. Let
$\eta_0>0$ such that
\beq\label{M1-3}
|f(x)-f(x')|<\dt/4\tforal f\in {\cal G},
\eneq
provided that ${\rm dist}(x, x')<\eta_0.$

Choose $\eta_0>\eta_1>0$ such that $B(x_i, \eta_1)$ intersects
with one and only one $S_i$ among $\{S_1,S_2,...,S_n\}.$

 Choose, for each $i,$ a  non-zero function $f_i\in C(X)$ with $0\le f\le 1$
whose support is in $B(x_i, \eta_1/2).$ Put
$$
d_0=\min\{ \inf\{\gamma(f_i)(\tau): \tau\in T(A)\}: 1\le i\le N\}.
$$
So $d_0>0.$
Put $\dt_1=\min\{\dt/8, \dt_0/4\}$ and put ${\cal G}_1={\cal
G}\cup\{1_{C(X)}\}\cup \{f_i: i=1,2,...,N\}.$

Now applying \ref{ML0}. We obtain a projection $p\in A$ and a
unital \hm\, $h_0: C(X)\to pAp$ such that
\beq\label{M1-4}
|\tau\circ h_0(g)-\gamma(g)(\tau)|&<&\dt_1\tforal g\in {\cal
G}_1\andeqn\\\label{M1-5}
 \tau\circ
h_0(\chi_{S_i})&<&\gamma(\chi_{S_i})(\tau)
\eneq
for all $\tau\in T(A),$ $i=1,2,...,n.$ Since $1_{C(X)}\in {\cal
G}_1,$ by (\ref{M1-4}),
\beq\label{M1-4+}
\tau(1-p)<\dt_1<\dt_0/4\tforal \tau\in T(A).
\eneq

Write $h_0(f)=\sum_{j=1}^Lf(\xi_j)q_j$ for all $f\in C(X),$ where
$\xi_j\in X$ and $\{q_1,q_2,...,q_L\}$ is a set of mutually
orthogonal projections with $\sum_{j=1}^Lq_j=p.$

Define
$$
e_i'=\sum_{\xi_j\in B(x_i, \eta_1/2)}q_j,\,\,\,i=1,2,...,N.
$$

 It follows from (\ref{M1-4}) that, for each
$i,$
\beq\label{Ml-6}
\tau(e_i')&\ge & \tau\circ h_0(f_i)\\
&>& \gamma(f_i)(\tau)-\dt_1> 3\dt_0/4\ge \tau(p)
\eneq
for all $\tau\in T(A).$ It follows that
$$
[e_i']\ge [p],\,\,\,i=1,2,...,N.
$$
There are projections $e_i\le e_i'$ such that
\beq\label{M1-6+}
[e_i]=[e_1]\ge [1-p],\,\,\,i=1,2,...,N.
\eneq

Define
\beq\label{M1-7}
h_1(f)&=&\sum_{\xi_j\not\in \cup_{i=1}^N B(x_i, \eta_1/2)}
f(\xi_j)q_j+\sum_{i=j}^N f(x_j)(e_i'-e_i)\andeqn\\
h(f)&=&\sum_{i=1}^N f(x_i)e_i\oplus h_1(f)
\eneq
for all $f\in C(X).$ Since $B(x_j, \eta_1/2)$ lies in one of
$S_i,$
$$
\tau\circ h(\chi_{S_i})=\tau\circ h_0(\chi_{S_i})\tforal \tau\in
T(A),
$$
$i=1,2,...,n.$ It follows from  (\ref{M1-5}) that (\ref{M1-1+})
holds. By the choice of $\eta_0,$ we also have
\beq\label{M1-8}
\|h_0(g)-h(g)\|<\dt/2\tforal h\in {\cal G}.
\eneq
Thus,  by (\ref{M1-4}),  (\ref{M1-1}) also holds.

\end{proof}

\begin{lem}\label{ML2}
Let $X$ be a compact metric space and let $A$ be a unital simple
\CA\, with tracial rank zero. Suppose that $\gamma: C(X)_{s.a}\to
{\rm Aff}(T(A))$ is a unital strictly positive linear map which is
compatible with a strictly positive \hm\, $\kappa_0: K_0(C(X))\to
K_0(A).$ Fix $\dt>0,$ $\eta>0,$ a finite subset ${\cal F}\subset
C(X)_{s.a.},$ an integer $N\ge 1,$ an $\eta$-dense subset
$\{x_1,x_2,...,x_N\}$ of $X,$ a finitely many mutually disjoint
clopen subset $S_1, S_2, ...,S_n\subset X$ with $\cup_{i=1}^n
S_i=X,$ a finite subset set $\{a_1,a_2,...,a_n\}\subset A$ of
mutually orthogonal projections  with
$$
0<a_i< \kappa_0([\chi_{S_i})]),\,\,\,i=1,2,...,n,
$$
a finitely many mutually disjoint clopen subsets
$\{F_1,F_2,...,F_{n_1}\}$ of $X$ with $\cup_{i=1}^{n_1} F_i=X,$
and a projection $p$ with $\tau(p)=\tau(\sum_{i=1}^na_i)$ for all
$\tau\in T(A).$

There is a projection $q\in A$ such that $[p]\le [q]$ and a unital
\hm\, $h: C(X)\to qAq$ with finite dimensional range such that
\beq\label{M2-1}
|\tau\circ h(g)-\gamma(g)(\tau)|&<&\dt\tforal g\in {\cal F}\tand
\tau\in T(A), \tand\\\label{M2-1+} \tau\circ
h(\chi_{F_i})&<&\gamma(\chi_{F_i})(\tau)\tforal \tau\in T(A),
\eneq
$i=1,2,...,n,$
\beq\label{M2-2}
h(f)=\sum_{i=1}^N f(x_i)e_i\oplus h_1(f)\tforal f\in C(X),
\eneq
where $h_1: C(X)\to (1-\sum_{i=1}^Ne_i)A(1-\sum_{i=1}^N e_i)$ is a
unital \hm\, with finite dimensional range and
$\{e_1,e_2,...,e_N\}$ is a set of mutually orthogonal projections
such that $[e_i]=[e_1]\ge [1-p],$ $i=1,2,...,N.$

Moreover, there exists a projection $p'\in q$ such that
\beq\label{M2-3}
p'h(f)&=&h(f)p'\tforal f\in C(X)\andeqn
\eneq
\beq\label{M2-4}
[h(\chi_{S_j})p']&=&[a_j],\,\,\,j=1,2,...,n.
\eneq
\end{lem}

\begin{proof}
Let
$$
d_0=\min\{ \inf\{\tau(\kappa_0([\chi_{S_i})])-[a_i]): \tau\in
T(A)\}: 1\le i\le n\}
$$
and let
$$
d_1=\inf\{ \tau(1-p): \tau\in T(A)\}.
$$
Then $d_0, d_1>0.$
Define $\dt_1=\min\{\dt/4, d_0/2,d_1/2\}$ and ${\cal G}_1={\cal
F}\cup \{1_{C(X)}, \chi_{S_i}, i=1,2,...,n\}.$ By applying
\ref{ML1}, we obtain a projection $q\in A$ and a unital \hm\, $h:
C(X)\to qAq$ with finite dimensional range satisfying the
following:
\beq\label{M2-5}
|\tau\circ h(g)-\gamma(g)(\tau)|&<&\dt_1\tforal g\in {\cal G}_1,\\\label{M2-5+}
\tau\circ h(\chi_{F_j})&<&\gamma(\chi_{F_j})(\tau)
\eneq
for all $\tau\in T(A),$ $j=1,2,...,n_1,$
\beq\label{M2-6}
h(f)=\sum_{k=1}^Nf(x_k)e_i\oplus h_1(f)\tforal f\in C(X),
\eneq
where $\{e_1,e_2,...,e_N\}$ is a set of mutually orthogonal and
mutually equivalent projections such that $[e_1]\ge [1-q],$ and where $h_1: C(X)\to (q-\sum_{k=1}^Ne_k)A(q-\sum_{k=1}^Ne_k)$ is a unital \hm\, with finite dimensional range.

Since $1_{C(X)}\in {\cal G}_1$, by the choice of $\dt_1,$  we
conclude that $[p]\le [q].$

Moreover, by (\ref{M2-5}),
\beq\label{M2-7}
\tau\circ h(\chi_{S_i})>\tau(a_i)\tforal \tau\in
T(A),\,\,\,i=1,2,...,n.
\eneq
Write
$$
h(f)=\sum_{s=1}^L f(\xi_s)E_s\tforal f\in C(X),
$$
where $\xi_s\in X$ and $\{E_1, E_2,...,E_L\}$ is a set of mutually
orthogonal projections such that $\sum_{s=1}^L E_s=q.$ By
(\ref{M2-7}), one has
\beq\label{M2-8}
\sum_{\xi_s\in S_i} E_s\ge a_i,\,\,\, i=1,2,...,n.
\eneq
For each $i,$ by the Riesz Interpolation Property (\cite{Z1}), there is a
projection $E_s'\le E_s$ for which $x_s\in S_i$ such that
\beq\label{M2-9}
[\sum_{\xi_s\in S_i}E_s']=[a_i].
\eneq

Put $p'=\sum_{s=1}^L E_s'$ Then
\beq
p'h(f)&=&h(f)p'\tforal f\in C(X)\andeqn
\eneq
\beq
[h(\chi_{S_i})p']=[a_i],\,\,\,i=1,2,...,n.
\eneq

\end{proof}

\begin{thm}\label{TKL1}
Let $X$ be a compact subset of a finite CW complex and let $A$ be
a unital simple \CA\, with $TR(A)=0.$  Suppose that $\kappa\in
KL_e(C(X), A)^{++}$ and suppose that there exists a unital
strictly positive linear map $\gamma: C(X)_{s.a}\to {\rm
Aff}(T(A))$ which is compatible with $\kappa.$
 Then there exists a
unital monomorphism $\phi: C(X)\to A$ such that
$$
[\phi]=\kappa\,\,\,\text{in}\,\,\, KL(C, A).
$$
\end{thm}

\begin{proof}


Suppose that $X\subset Y,$ where  $Y$ is a finite CW complex.
Let $X_n\subset Y$ be a decreasing sequence of finite CW complexes
for which \ref{D1} holds. Suppose that
$p_{n,1},p_{n,2},...,p_{n,r(n)}$ are  mutually orthogonal
projections of $C(X_n)$ which correspond to the connected
components of $X_n.$ It is clear that we may assume that each
connected component of $X_n$ contains at least one point of $X.$
This implies that $[s_n]\in KK(C(X_n), C(X))^{++}.$ It follows
that
\beq\label{TKL1-2}
[s_n]\times \kappa\in KK(C(X_n), A)^{++}.
\eneq

Let $\{{\cal F}_n\}$ be an increasing sequence of finite subsets
of $C(X)$ whose union is dense in $C(X).$ Let $\{\eta_n\}$ be a
decreasing sequence of positive numbers with
$\lim_{n\to\infty}\eta_n=0,$ $\{{\cal P}_n\}$ be an increasing
sequence of finite subsets of $\underline{K}(C)$ whose union is
$\underline{K}(C),$ let $\{k(n)\},\, \{N(n)\}\subset \N$ be two
sequences of integers such that $k(n), N(n)\nearrow \infty,$ and
$\{x_{n,1},x_{(n,2},...,x_{n,N(n)}\}$ be $\eta_n$-dense
 subsets of $X$ which satisfy the requirements of \ref{KLL1} and \ref{RM1}
for corresponding $\ep_n=1/2^{n+2}$ and  ${\cal F}_n.$

By passing to a subsequence if necessary, we may assume that there
is a finite subset ${\cal F}_n'\subset C(X_{k(n+1)})$ such that
$s_{k(n+1)}({\cal F}_n')={\cal F}_n$ and a finite subset ${\cal
Q}_{k(n)}\subset \underline{K}(C(X_{k(n)}))$ such that
$[s_{k(n)}]({\cal Q}_{k(n)})={\cal P}_n,$ $n=1,2,....$ We may
assume that $1_{C(X_{k(n)})}\in {\cal F}_n',$ without loss of
generality.

Set $\kappa_n=[s_{k(n)}]\times \kappa.$  Note that
$\kappa_n([1_{C(X_{k(n)})}])=[1_A].$


Let $\dt_n$ (in place of $\dt$), ${\cal G}_n'\subset C(X)_{s.a}$
( in place of ${\cal G}$), $S_{1,n}, S_{2,n},...,S_{m(n),n}$ (in
place of $\{S_1, S_2,...,\}$) be a set of disjoint clopen
subsets of $X$  with $\cup_{i=1}^{m(n)}S_i=X$ required by
\ref{Lcd} for $\ep_n$ and ${\cal F}_n,$ $n=1,2,....$ We may assume
that $1_{C(X)}\in {\cal G}_n',$ $n=1,2,....$

By taking a refinement of the clopen partition of $X,$ we may
assume that $s_n(p_{n,i})$ is a finite sum of functions in
$\{\chi_{S_{j,n}}: 1\le j\le m(n)\},$ $i=1,2,...,r(n).$

Let ${\cal G}_n\subset C(X_{k(n)})_{s.a.}$ be a finite subsets for
which $s_{k(n)}({\cal G}_n)={\cal G}_n',$ $n=1,2,....$

By applying \ref{ML1}, we obtain a projection $P_1\in A$ and a
unital \hm\, $\Phi_1': C(X)\to P_1AP_1$ such that
\beq\label{T1-10}
|\tau\circ \Phi_1'(g)-\gamma(g)(\tau)|&<&\dt_1/2\tforal g\in {\cal
G}_1',\\
\tau\circ \Phi_1'(\chi_{S_{j,1}})&<&\gamma(\chi_{S_{j,1}})(\tau)
\eneq
for all $\tau\in T(A),$ $i=1,2,...,m(1),$ and
\beq\label{T1-11}
\Phi_1'(f)=\sum_{i=1}^{N(1)} f(x_{1,i})e_i^{(1)}\oplus
\Phi_{0,1}'(f)\tforal f\in C(X),
\eneq
where $\{e_1^{(1)},e_2^{(1)},...,e_{N(1)}^{(1)}\}$ is a set of
mutually orthogonal and mutually equivalent projections with
$[e_1]\ge [(1-P_1)]$ and where $ \Phi_{0,1}':
C(X)\to(P_1-\sum_{i=1}^{N(1)} e_i^{(1)})A((P_1-\sum_{i=1}^{N(1)}
e_i^{(1)})$ is a unital \hm\, with finite dimensional range.  Note
also, since $1_{C(X)}\in {\cal G}_1',$ $\tau(1-P_1)<\dt_1/2$ for
all $\tau\in T(A).$

It follows from \ref{KLL2} that there is a unital monomorphism
$\phi_1': C(X_{k(1)})\to (1-P_1)A(1-P_1)$ such that
\beq\label{T1-12}
[\phi_1']+[\Phi_1'\circ s_1]=\kappa_1\,\,\,{\rm in}\,\,\,
KK(C(X_{k(1)}), A).
\eneq

Define $\phi_1=\phi_1'+\Phi_1'\circ s_1.$

Suppose that, for $1\le m\le n,$ there are unital \hm s  $\phi_m': C(X_{k(m)})\to (1-P_m)A(1-P_m)$  and
$\Phi_m': C(X)\to P_mAP_m$  and a unital (injective) \hm\, 
$\phi_m=\phi_m'+\Phi_m'\circ s_{k(m)}$ such that
\begin{enumerate}
\item there are mutually orthogonal and mutually equivalent
projections $e_1^{(m)}, e_2^{(m)},...,e_{N(m)}^{(m)}\in P_mAP_m$
for which $[e_1^{(m)}]\ge [1-P_m],$ and
$$
\Phi_m'(f)=\sum_{i=1}^{N(m)}f(x_{m,i})e_i^{(m)}\oplus
\Phi_m^{(0)}(f) \tforal f\in C(X)
$$
where  $\Phi_m^{(0)}: C(X)\to
(P_m-\sum_{i=1}^{N(m)}e_i^{(m)})A(P_m-\sum_{i=1}^{N(m)}e_i^{(m)})$
is a unital \hm\, with finite dimensional range;

\item  $\tau\circ
\Phi'(\chi_{S_{j,m}})<\gamma(\chi_{S_{j,m}})(\tau)\tforal \tau\in
T(A),$ $j=1,2,...,m(m);$

\item $|\tau\circ \Phi_m'(g)-\gamma(g)(\tau)|<\dt_m/2\tforal g\in
{\cal G}_m'$ and for all $\tau\in T(A);$

\item $[P_{m+1}]\ge [P_m]$ in $K_0(A)$ and $\tau(1-P_m)<\dt_m/2
\tforal \tau\in T(A);$

\item there is a projection $P_{m+1}'\le P_{m+1}$ such that
$P_{m+1}'\Phi_{m+1}=\Phi_{m+1}'P_{m+1}$ and
$$[\Phi_{m+1}'(\chi_{S_{j,m}})P_{m+1}']=[\Phi_m'(\chi_{S_{j,m}})]\,\,\,
{\rm in}\,\,\,K_0(A),j=1,2,...,m(m);$$

 \item $\phi_m'$ is a unital monomorphism;

\item $[\phi_m]=[\phi_m']+[\Phi_m'\circ s_{k(m)}]=\kappa_m;$

 \item there exists a unitary $u_m\in A$ such that
$$
{\rm ad}\, u_m\circ \phi_{m+1}\circ
s_{k(m), k(m+1)}\approx_{1/2^{m+1}} \phi_{m}\,\,\,{\rm on}\,\,\, s_{k(m)}^{-1}({\cal F}_m
),\,\,\,m=1,2,...,n-1.
$$

\end{enumerate}

It follows from \ref{ML2} that there is a projection $P_{n+1}\in
A$ and a unital \hm\, $\Phi_{n+1}': C(X)\to P_{n+1}AP_{n+1}$
satisfying the following:

\begin{enumerate}
\item there are mutually orthogonal and mutually equivalent
projections $e_1^{(n+1)}, e_2^{(n+1)},...,e_{N(n+1)}^{(n+1)}\in
P_{n+1}AP_{n+1}$ for which $[e_1^{(n+1)}]\ge [1-P_{n+1}],$ and
$$
\Phi_{n+1}'(f)=\sum_{i=1}^{N(n+1)}f(x_{n+1,i})e_i^{(n+1)}\oplus
\Phi_{n+1}^{(0)}(f) \tforal f\in C(X)
$$
where  $\Phi_{n+1}^{(0)}: C(X)\to
(P_{n+1}-\sum_{i=1}^{N(n+1)}e_i^{(n+1)})A(P_{n+1}-\sum_{i=1}^{N(n+1)}e_i^{(n+1)})$
is a unital \hm\, with finite dimensional range;
\item  $\tau\circ
\Phi_{n+1}'(\chi_{S_{j,n+1}})<\gamma(\chi_{S_{j,n+1}})(\tau)\tforal
\tau\in T(A),$ $j=1,2,...,m(n+1);$

\item $|\tau\circ
\Phi_{n+1}'(g)-\gamma(g)(\tau)|<\dt_{n+1}/2\tforal g\in {\cal
G}_{n+1}'$ and for all $\tau\in T(A);$

\item $[P_{n+1}]\ge [P_n]$ in $K_0(A)$ and
$\tau(1-P_{n+1})<\dt_{n+1}/2 \tforal \tau\in T(A);$

\item there is a projection $P_{n+1}'\le P_{n+1}$ such that
$P_{n+1}'\Phi_{n+1}'=\Phi_{n+1}' P_{n+1}$ and
$$[\Phi'_{n+1}(\chi_{S_{n, j}})P_{n+1}']=[\Phi'_n(\chi_{S_{n,
j}})]\,\,\,{\rm in}\,\,\,K_0(A),j=1,2,...,m(n).$$
\end{enumerate}

It follows from \ref{KLL2} that there is a unital
monomorphism\linebreak  $\phi_{n+1}': C(X_{k(n+1)})\to
(1-P_{n+1})A(1-P_{n+1})$ such that
\beq\label{T1-13}
[\phi_{n+1}']=\kappa_{n+1}-[\Phi_{n+1}'\circ s_{k(n+1)}]\,\,\,{\rm
in}\,\,\, KK(C(X_{k(n+1)},A)
\eneq
Define $\phi_{n+1}=\phi_{n+1}'+\Phi_{n+1}'\circ s_{k(n+1)}.$

Thus $\phi_{n+1}',$ $\phi_{n+1}'$ and $\phi_{n+1}$ satisfy (1),
(2), (3), (4), (5), (6) and (7).

To complete the induction, define  $\Phi_{n+1}'': C(X)\to
P_{n+1}'AP_{n+1}'$ by
$\Phi_{n+1}''(f)=P_{n+1}'\Phi_{n+1}'(f)P_{n+1}'$ for all $f\in
C(X).$ By (3) and (4),
$$
|\tau\circ \Phi_{n+1}''(g)-\gamma(g)(\tau)|<\dt_{n+1}/2\tforal g\in
{\cal G}_n
$$
for all $\tau\in T(A).$ Note that, by (5), $[P_{n+1}']=[P_n].$
There is a unitary $w_n\in U(A)$ such that
$$
w_n^*P_{n+1}'w_n=P_n.
$$
Thus, by (5) and (3), and by applying \ref{Lcd}, there exists a
unitary $v_n\in U(P_nAP_n)$ such that
\beq\label{T1-14}
{\rm ad}\, v_n\circ {\rm ad}\, w_n\circ
\Phi_{n+1}''\approx_{\ep_n} \Phi_n'\,\,\,{\rm
 on}\,\,\, {\cal F}_n.
 \eneq
 Denote $\Psi_{n+1}'=P'_{n+1}\Phi_{n+1}P'_{n+1}$ and $\Psi_{n+1}={\rm ad}\, w_n\circ \Psi_{n+1}'.$ 
 Let $\phi_{n+1}''={\rm ad}\, w_n\circ \phi_{n+1}'\oplus \Psi_{n+1}.$ 
 Now consider $\phi_n'$ and  $\phi_{n+1}''\circ s_{k(n), k(n+1)}.$  By (7) 
and (\ref{T1-13}),
we have
$$
[\phi_{n+1}''\circ s_{k(n), k(n+1)}]|_{{\cal Q}_{k(n)}}=[\phi_n']|_{{\cal Q}_{k(n)}}.
$$
It follows from (1) and 
\ref{KLL1} that there exists a unitary $V_n\in U(A)$ such that
\beq\label{T1-15}
{\rm ad}\, V_n\circ (\phi_{n+1}''\circ s_{k(n),
k(n+1)}\oplus \Phi_n' \circ s_{k(n)})\approx_{\ep_n} \phi_n'\oplus
\Phi_n'\circ s_{k(n)}\,\,\,{\rm on}\,\,\, s_{k(n)}^{-1}({\cal
F}_n).
\eneq
Define $u_n=w_n (v_n+(1-P_n)) V_n.$ Then, by (\ref{T1-14}) and
(\ref{T1-15}),
\beq\label{T1-16}
{\rm ad}\,u_n\circ \phi_{n+1}\approx_{2\ep_n}\phi_n\,\,\,{\rm
on}\,\,\, s_{k(n)}^{-1}({\cal F}_n).
\eneq
Note $2\ep_n=1/2^{n+1}.$

This concludes the induction.

Define $\psi_1=\phi_1$ and $\psi_{n+1}={\rm ad}\,u_n\circ
\phi_{n+1},$ $n=1,2,....$ Then, by (8) above,
\beq\label{T1-5}
\|\psi_{n}(c)-\psi_{n+1}\circ
s_{k(n),k(n+1)}(c)\|<1/2^{n+2}\tforal c\in s_{k(n)}^{-1}({\cal
F}_n),
\eneq
$n=1,2,....$


Fix $m$ and $f\in {\cal F}_m,$ let $g\in s_{k(m)}^{-1}({\cal
F}_m)$ such that $s_{k(m)}(g)=f.$

It follows that $\{\psi_n\circ s_{m,n}(g)\}_{n\ge m}$ is a Cauchy
sequence by (\ref{T1-5}).

Note that if $g'\in s_{k(m)}^{-1}({\cal F}_m)$ such that
$s_{k(m)}(g)=s_{k(m)}(g'),$ then, for any $\ep>0,$  there exists
$n\ge m$ such that
$$
\|s_{k(m), k(n)}(g)-s_{k(m), k(n)}(g')\|<\ep.
$$

Thus  $h(f)=\lim_{n\to\infty}\psi_n\circ s_{m,n}(g)$ is
well-defined. It is then easy to verify  that $h$ defines a unital
\hm\, from $C(X)$ into $A.$ Since each $\phi_n$ is injective, it
is easy to check that $h$ is also injective.

If  $x\in {\cal Q}_m,$ then by (7) above,
$$
 [h]\circ [s_{k(m)}](x)=\kappa_n\circ [s_{k(m),
k(n)}](x)=\kappa\circ [s_{k(m)}](x).
$$
Therefore
$$
[h]=\kappa\,\,\,{\rm in}\,\,\, KL(C,A).
$$

It is also easy to check from (3) and (4) that
\beq
\tau\circ h(g)=\gamma(\check{g})(\tau)\tforal g\in C(X)_{s.a}
\eneq
and for all $\tau\in T(A).$

\end{proof}

\section{AH-algebras}

\begin{lem}\label{2L}
Let $X$ be a compact subset of a finite CW complex, let
$C=PM_k(C(X))P,$ where $P\in M_k(C(X))$ is a projection, and let
$A$ be a unital simple \CA\, with tracial rank zero. Suppose that
$\kappa\in KL_e(C, A)^{++}$ and suppose that $\gamma: C_{s.a.}\to
{\rm Aff}(T(A))$ is  a unital positive linear map which is
compatible with $\kappa.$ Then, for any $\ep>0$ and finite subset
 ${\cal F}\subset  C(X),$ there is a unital monomorphism
 $\phi: C\to A$ such that
 \beq\label{KLL4}
 [\phi]=\kappa\andeqn \tau\circ \phi(f)=\gamma(f)(\tau)\tforal \tau\in T(A).
 \eneq

\end{lem}

\begin{proof}
It is clear that the case that $C=M_k(C(X))$ follows from
\ref{TKL1} immediately. For the general case, there is an integer
$d\ge 1$ and a projection $p\in M_d(C)$ such that $pM_d(C)p\cong
M_m(C(X))$ for some integer $m\ge 1.$ Thus the general case is
reduced to the case that $C=M_m(C(X)).$

\end{proof}

\begin{thm}\label{T2}
Let $C$ be a unital AH-algebra and let $A$ be a unital simple
\CA\, with $TR(A)=0.$ Suppose that $\kappa\in \mathrm{KL}_e(C,
A)^{++}.$  Suppose also that there is a unital strictly positive
linear map $\gamma: C_{s,a}\to {\rm Aff}(T(A))$ which is
compatible with $\kappa.$ Then there is a monomorphism $\phi: C\to
A$ such that
\beq\label{KL-1}
[\phi]&=&\kappa\,\,\,\textrm{in}\,\,\,
\mathrm{KL}(C,A)\andeqn\\
\tau\circ \phi(c)&=&\gamma(c)(\tau)
\eneq
for all $c\in C_{s.a.}$ and $\tau\in T(A).$
\end{thm}

\begin{proof}

We may write $C=\overline{\cup_{n=1}^{\infty}C_n},$ where $C_n
=P_nM_k(C(X_n))P_n,$ where $X_n$ is a compact subset of a finite
CW complex and $P_n\in M_k(C(X_n))$ is a projection. We may also
assume that $1_{C_n}=1_C$ for all $n.$ Denote by $\imath_n: C_n\to
C$  the embedding, $n=1,2,....$

Define
$$
\kappa_n=\kappa\circ [\imath_n]\andeqn \gamma_n=\gamma\circ
(\imath_n)_{\sharp}
$$
$n=1,2,....$ Since $\imath_n$ is injective $\kappa_n\in KL_e(C_n,
A)^{++}$  and $\gamma_n$ is unital strictly positive. It is also
clear that $\gamma_n$ is compatible with $\kappa_n,$ since $\gamma$
is compatible with $\kappa.$ It follows from \ref{2L} that there
is a sequence of unital monomorphisms $\phi_n: C_n\to A$ such that
\beq\label{T2-3}
[\phi_n]=\kappa_n\andeqn \tau\circ \phi_n(c)=\gamma_n(c)(\tau)
\eneq
for all $c\in C_{s.a.}$ and $\tau\in T(A).$

Let $\{{\cal F}_n\}$ be an increasing sequence of finite subsets
of $C$ whose union is dense in $C.$ By passing to a subsequence,
if necessary, without loss of generality, we may assume that
${\cal F}_n\subset C_n.$

It follows (from  2.3.13 of \cite{Lnbk}, for example)  that there
is, for each $n,$ a unital completely positive linear map $L_n:
C\to A$ such that
\beq\label{T2-4}
L_n\approx_{1/2^{n+1}} \phi_n\circ \imath_n\,\,\,{\rm
on}\,\,\,{\cal F}_n.
\eneq

It follows from Lemma \ref{2L}, by passing to a subsequence again and by
applying (\ref{T2-3}), there is a sequence of unitaries $u_n$ and a
subsequence of $\{k(n)\}$ such that
\beq\label{T2-5}
{\rm ad}\, u_n\circ L_{k(n+1)}\approx_{1/2^n} L_{k(n)}\,\,\,{\rm
on}\,\,\, {\cal F}_n,
\eneq
$n=1,2,....$ Define $\psi_1=L_1,$ $\psi_{n+1}={\rm ad}\, u_n\circ
L_{n+1},$ $n=1,2,....$  Note that $\{\psi_n(c)\}$ is a Cauchy
sequence in $A$ for each $c\in {\cal F}_m.$ Define
$h(c)=\lim_{n\to\infty}\psi_n(c).$ It is easy to see that $h$
gives a unital \hm\, from $C$ into $A.$ Moreover, for each $x\in
\cup_{n=1}^{\infty}{\cal F}_n,$
\beq\label{T2-6}
h(x)=\lim_{n\to\infty} {\rm ad}\, u_n \circ \phi_{k(n)}\circ
\imath_{k(n)}\circ\cdots \circ \imath_n(x).
\eneq
Since each $\phi_n$ is injective, it follows that $h$ is a
monomorphism. From (\ref{T2-6}) and (\ref{T2-3}), we have
$$
[h]=\kappa\,\,\,{\rm as\,\,\, well \,\,\, as}\,\,\, \tau\circ
h(c)=\gamma(c)(\tau)
$$
for all $c\in C_{s.a.}$ and $\tau\in T(A).$

\end{proof}

\begin{cor}\label{CCX}
Let $X$ be a compact metric space and let $A$ be a unital simple
\CA\, with tracial rank zero. Suppose that $\kappa\in
\mathrm{KL}_e(C(X), A)^{++}.$ Suppose also that there is a unital
strictly positive linear map $\gamma: C_{s,a}\to {\rm Aff}(T(A))$
which is compatible with $\kappa.$ Then there is a monomorphism
$\alpha: C\to A$ such that
\beq\label{CKL-1}
[\alpha]&=&\kappa\,\,\,\textrm{in}\,\,\,
\mathrm{KL}(C(X),A)\andeqn\\
\tau\circ \phi(c)&=&\gamma(c)(\tau)
\eneq
for all $c\in C(X)_{s.a.}$ and $\tau\in T(A).$
\end{cor}

\begin{exm}\label{Ex}
{\rm Let $X=\{{-1\over{n}}: n\in \N\}\cup [0,1]\cup
\{1+{1\over{n}}: n\in \N\}\subset [-1, 2].$ Put $C=C(X).$ Then
$$
K_0(C(X))=C(X, \Z).
$$
Take two sequences of positive rational numbers $\{a_n\}$ and
$\{b_n\}$ such that $\sum_{n=1}^{\infty} a_n =1-\sqrt{2}/2$ and
$\sum_{n=1}^{\infty}=\sqrt{2}/2.$

Define a unital positive linear functional $F: C(X)\to \R$ as
follows:
$$
F(f)=\sum_{n\in\N} a_nf({-1\over{n}})+\sum_{n\in \N} b_n
f({1\over{n}})\tforal f\in C(X).
$$
Let $D_0=F(C(X, \Z)).$ Note that, if $S$ is a clopen subset which
does not contain $[0,1],$ then $F(S)\in \Q.$ If $S\supset [0,1], $
Then
$$
F(S)=1-F(S_1)
$$
for some clopen subset $S_1\subset X$ which does not intersect
with $[0,1].$ It follows that $D_0\subset \Q.$

This gives a unital positive linear map $F_*: C(X, \Z)\to \Q.$ Let
$p\in C(X)$ be a projection whose support $\Omega$ has a non-empty
intersection with $[0,1].$ Since $\Omega$ is clopen,
$\Omega\supset [0,1].$ It follows that there is $N\ge 1$ such that
${1\over{k}}\in \Omega$ for $|k|\ge N.$ It follows that
$$
F(p)\ge \sum_{|k|\ge N} {1\over{2^{|k|+1}}}>0.
$$
From this one sees that $F_*$ is strictly positive.

Let $A$ be a unital simple AF-algebra with
$$
(K_0(A), K_0(A), [1_A])=(\Q, \Q_+,1).
$$
There is an element $\kappa\in KL(C(X), A)$ such that
$$
\kappa|_{K_0(C(X))}=F_*.
$$
Thus $\kappa(K_0(C(X))_+\setminus \{0\})\subset K_0(A)_+\setminus
\{0\}.$ In other words,  $\kappa\in KL_e(C,A)^{++}.$

Suppose that $\gamma: C_{s.a}\to {\rm Aff}(T(A))=\R$ is unital and
positive such that
$$
\gamma(\check{p})(\tau)=\tau(\kappa([p]))
$$
for all projections $p\in C$ and $\tau\in T(A).$ Consider a
positive continuous function $f\in C(X)$ with $0\le f\le 1$ whose
support $S$ is an open subset of $(0,1).$ Consider
projection $p_n(t)=0$ if $t\not\in [-1/n, 1+1/n]\cap X$ and
$p_n(t)=1$ if $t\in [-1/n, 1+1/n]\cap X.$ Then
$$
f\le p_n,\,\,\, n=1,2,....
$$
It follows that, for all $\tau\in T(A),$
\beq
\gamma (\check{f})(\tau)&\le &\gamma(\check{p}_n)(\tau)\\
&<&\sum_{|k|\ge n} (a_k+b_k)\to 0
\eneq
as $|n|\to \infty.$ It follows that
$$
\gamma(\check{f})(\tau)=0\tforal \tau\in T(A).
$$
This shows that $\gamma$ is not strictly positive.

In particular, there is {\it no} unital monomorphism $\phi:
C(X)\to A$ such that $[\phi]=\kappa.$

How about \hm s?  Suppose that there exists a unital \hm\, $h:
C(X)\to A$ such that $[h]=\kappa.$ Let $f\in C(X)_+$ be so that
its support is contained in $[0,1].$ Then, as shown above,
$\tau(h(f))=0$ for $\tau\in T(A).$  Since $A$ is simple, this
implies that $h(f)=0.$ It is then easy to see that
$$
{\rm ker}\,h=\{f\in C(X): f|_{X\setminus (0,1)}=0\}.
$$
Thus $C/{\rm ker}\,h\cong C(Y),$ where $Y=X\setminus (0,1).$ Let
$\phi: C(Y)\to A$ be the unital \hm\, induced by $h.$ Then $\phi$
is a monomorphism.  Let
$$
Y_1=\{1+1/n: n\in \N\}\cup\{1\}\andeqn Y_2=\{-1/n: n\in
\N\}\cup\{0\}.
$$
Then $Y_1$ and $Y_2$ are  clopen subsets of $Y.$ Let $p_i$ be the
projection corresponding to $Y_i,$ $i=1,2.$ Then
$$
\tau(p_1)\ge \sum_{n=1}^{\infty} b_n=1-\sqrt{2}/2\andeqn\\
\tau(p_2)\ge \sum_{n=1}^{\infty} a_n=\sqrt{2}/2
$$
for $\tau\in T(A).$ Since $\tau(p_1)+\tau(p_2)=1,$ it follows that
$$
\tau(p_1)=1-\sqrt{2}/2\andeqn \tau(p_2)=\sqrt{2}/2.
$$
This is impossible since $K_0(A)=\Q.$

}
\end{exm}

From this we arrive at the following conclusion:

\begin{prop}\label{cex}
There are compact metric spaces $X$ with dimension one, unital
simple AF-algebras $A$ with unique tracial states and  $\kappa\in
KL_e(C,A)^{++}$ which has no  strictly positive affine map from
${\rm Aff}(T(C(X))$ to ${\rm Aff}(T(A))$ compatible with $\kappa.$

Furthermore, there is no unital \hm\, $\phi: C(X)\to A$ such that
$[\phi]=\kappa$ in $KL(C,A).$

\end{prop}




\

\begin{df}\label{DMT}
Let $C$ be a unital AH-algebra which admits a faithful tracial
state and let $A$ be a unital simple \CA\, with
$T(A)\not=\emptyset.$

Denote by $ KLT(C,A)^{++} $ the set of pairs $(\kappa, \lambda)$
where $\kappa\in KL(C,A)^{++}$  with $\kappa([1_C])=[1_A]$ and
$\lambda: T(A)\to T_{\mathtt{f}}(C)$ which is compatible with
$\kappa,$ i.e., $\lambda(\tau)(p)=\tau(\kappa([p])$ for all
projections $p\in M_{\infty}(C)$ and for all $\tau\in T(A).$

Denote by ${\rm Mon}_{au}^e(C,A)$ the set of approximately unitary
equivalent classes of unital monomorphisms from $C$ into $A.$

\end{df}

To conclude this note, combing the previous result in \ref{Lncd} (see \ref{Om}) and \ref{T2}, we state the following:

\begin{thm}\label{MT}
Let $C$ be a unital AH-algebra which admits a faithful tracial
state and let $A$ be a unital separable simple \CA\, with
$TR(A)=0.$ Then  map
$$
\Lambda: {\rm Mon}_{au}^e(C,A)\to  KLT(C,A)^{++}
$$
defined by $\phi\mapsto ([\phi], \phi_{T})$ is bijective.

\end{thm}

\end{document}